\documentclass[12pt]{amsart}
\usepackage{amssymb}
\usepackage{amssymb, pb-diagram}
\usepackage[all]{xy}

\newcommand{\Z}{{\mathbb{Z}}}

\numberwithin{equation}{section}

\newtheorem{lem}{Lemma}[section]
\newtheorem{corol}[lem]{Corollary}
\newtheorem{theor}[lem]{Theorem}
\newtheorem{prop}[lem]{Proposition}
\newtheorem{rema}[lem]{Remark}
\newtheorem{defi}[lem]{Definition}

\newtheorem{exem}[lem]{Example}

\begin{document}
\title[Algebras of quotients of path algebras]{Algebras of quotients of path algebras}

\author{M. Siles Molina}
\address{Departamento de \'Algebra, Geometr\'{\i}a y Topolog\'{\i}a, Universidad de M\'alaga, 29071 M\'alaga, Spain.}
\email{mercedes@agt.cie.uma.es}

\subjclass[2000]{Primary 16D70} \keywords{Leavitt path algebra,
path algebra, socle, algebra of quotients, Fountain-Gould right
order, Toeplitz algebra}

\begin{abstract} Leavitt path algebras are shown to be algebras
of right quotients of their corresponding path algebras. Using this
fact we obtain maximal algebras of right quotients from those
(Leavitt) path algebras whose associated graph satisfies that
every vertex connects to a line point (equivalently, the Leavitt
path algebra has  essential socle). We also introduce and
characterize the algebraic counterpart of Toeplitz algebras.
\end{abstract}

\maketitle

\section{Introduction and preliminaries}
Leavitt path algebras, natural generalizations of the algebras investigated by Leavitt in \cite{L}, are
algebraic versions of the Cuntz-Krieger algebras of directed graphs described in \cite{R} which, at the same
time, generalize Cuntz algebras (whose introduction and study was motivated by questions in physics). The
introduction of Leavitt path algebras in \cite{AA1} and \cite{AMP} has recently attracted the interest of a
significant number of algebraists as well as of analysts working on $C^*$-algebras. As a sample of this, let
us mention the notes of the ``Workshop on graph algebras'' (\cite{Notes}), held at the University of M\'alaga
(Spain), focused on both, the analytic and the algebraic part of graph algebras, through the history of the
subject and recent developments. Although the algebraic results look very similar (but are not exactly the
same, as shown in \cite{APS} and \cite{AP}), they require quite different techniques to be reached. We
recommend the reader the paper by Tomforde \cite{T}, where the author considers both areas of graph algebras.

 Roughly speaking, for a row-finite graph $E$ and a field $K$,  the
Leavitt path algebra $L_K(E)$ is the path $K$-algebra associated to $E$, modulo some relations (the so called
Cuntz-Krieger relations).

   One of the main interests of the researchers in Leavitt path algebras (as well as in graph C$^\ast$-algebras)
is to get a structure theory as deep as possible, and the developement of theories of algebras of quotients
allows to one to obtain a better understanding of this structure. This has been the motivation for the work
in this paper.

 Take, for example, a finite and acyclic graph $E$; then the Leavitt path algebra
$L_K(E)$ is semisimple and artinian, i.e, it is isomorphic to $\bigoplus_{i=1}^t \mathbb{M}_{n_i}(K)$ (see
\cite[Proposition 3.5]{AAS1}). Moreover, $L_K(E)$ is the maximal algebra of right quotients of the path
algebra $KE$ (\cite{E}). It is known that when the maximal algebra of right quotients of a semiprime algebra
$A$, let us call it $Q$, is semisimple and artinian, then $Q$ coincides with the classical algebra of right
quotients of $A$, hence it is natural to ask about when the path algebra $KE$ is semiprime. Another question
to be considered is if the Leavitt path algebra of an arbitrary graph is an algebra of right quotients of its
corresponding path algebra. A partial answer was given in \cite[Corollary 3.25]{E}: if the graph $E$ is
finite and acyclic, then $L_K(E)$ is an algebra of right quotients of the path algebra $KE$.

   In this paper we show that the path algebra $KE$ associated to
any graph $E$ is semiprime if and only if whenever there is a path joining two vertices, there exists another
one from the range to the source of the first one (Proposition \ref{SemipOfTheGraphAlgebra}). This implies
that for $E$ a finite and acyclic graph the Leavitt path algebra $L_K(E)$ is not the classical algebra of
right quotients of the path algebra $KE$ (except if there are no edges joining vertices), although it is the
classical algebra of right quotients of any semiprime subalgebra containing the path algebra (Proposition
\ref{dmaximalcasofinito}). Moreover, for a general graph $E$ we obtain that $L_K(E)$ is an algebra of right
quotients of the path algebra $KE$ (Proposition \ref{rqa}), result that extends that of \cite{E}.

There exists a notion of order in nonunital rings which was introduced by Fountain and Gould in \cite{FG1}
and extended by \'Ahn and M\'arki to one-sided orders (see \cite{AM1}). These authors developed, some years
later, a general theory of Fountain-Gould one-sided orders (see \cite{AM2} and the references therein). This
localization was carried out by considering generalized inverses like group inverses, and agrees with the
classical one (where usual inverses are taken) when the ring is semiprime and coincides with its socle
(\cite[Corollary 3.4]{GS}). If a ring $R$ has a Fountain-Gould right quotient ring $Q$, then $Q$ is unique up
to isomorphisms (see \cite[Theorem 5.9]{VG} or \cite[Corollary 3]{AM2}).

It was shown in \cite{AM2}, and later in \cite{GS}, that the maximal ring of right quotients provides an
appropriate framework where to settle these right orders, specially when the ring $R$ is a Fountain-Gould
right order in a semiprime ring coinciding with its socle. On the other hand, since the socle of a Leavitt
path algebra has been studied in \cite{AMMS}, we have the required tools to establish results on right
orders in Leavitt path algebras which satisfy the descending chain condition on principal one-sided ideals
(i.e., are semisimple).

   We prove that for an acyclic graph  $E$ any semiprime
subalgebra $A$ such that $KE\subseteq A \subseteq L_K(E)$ is a Fountain-Gould right order in the Leavitt path
algebra $L_K(E)$ (Proposition \ref{fountaingould}) and characterize these algebras $A$ when the hereditary
closure of the set of line points is the set of all vertices, that is, when $L_K(E)$ coincides with its socle
(Theorem \ref{leftlocalgoldie}).

In Section 4, we obtain that the socle of a Leavitt path algebra is an essential ideal if and only if every
vertex connects to a line point (Theorem \ref{socleessential}). This allows to reduce the study of maximal
algebras of quotients of the path algebra $KE$ (for $E$ as before) to the maximal algebras of quotients of
acyclic graphs $E^\prime$ such that $L_K(E^\prime)$ is semisimple, i.e., coincides with its socle (Theorem
\ref{elmaximaleseldelzocalo}). The maximal symmetric algebra of quotients and the maximal algebra of right
quotients of path algebras of locally finite graphs without cycles and such that the set of sinks is a
maximal antichain have been described in \cite[Chapter 3]{E}.

Finally, we define the algebraic counterpart of the Toeplitz algebra, as the Leavitt path algebra $T$ whose
graph is the following:
\medskip
$$
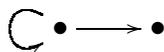
\xymatrix{ {\bullet} \ar@(ul,dl) \ar[r]&{\bullet}}$$
\medskip
Similarly to what happens in the analytic case, $T$ has an
essential ideal, which is the socle of $T$ in our setting, and
 there exists an exact sequence
$$0 \rightarrow Soc(L_K(E))\rightarrow T\rightarrow K[x, x^{-1}]\rightarrow 0.$$

Moreover, the Toeplitz algebra $T$ is sandwiched as follows:
$$\mathcal{M}_\infty(K)\subseteq T \subseteq \mathrm{RCFM}(K),$$
where $\mathcal{M}_\infty(K)$ denotes the algebra of matrices of infinite size with only a finite number of
nonzero entries.

As shown in the analytic context, $T$ has not a unique representation. It can be described as the Leavitt
path algebra $L_K(E)$, for $E$ a graph  such that $E^0=F^0\cup \{v\}$, $E^1=F^1\cup\{e, e_1, \dots, e_n\}$,
where $e$ has range and source $v$, $s(e_i)=v$, $r(e_i)\in F^0$, and for every $f\in F^0$ its range and
source are the corresponding as edges of the graph $F$, being $F$ an acyclic graph such that every vertex
connects to a line point. All of these results are collected in Theorem \ref{descriptionoftoeplitzalgebras}.

We start the preliminaries by recalling the definitions of path algebra and Leavitt path algebra. A
\emph{(directed) graph} $E=(E^0,E^1,r,s)$ consists of two countable sets $E^0,E^1$ and maps $r,s:E^1 \to
E^0$. The elements of $E^0$ are called \emph{vertices} and the elements of $E^1$ \emph{edges}. If $s^{-1}(v)$
is a finite set for every $v\in E^0$, then the graph is called \emph{row-finite}. Throughout this paper we
will be concerned only with row-finite graphs. If $E^0$ is finite then, by the row-finite hypothesis, $E^1$
must necessarily be finite as well; in this case we simply say that $E$ is \emph{finite}. A vertex which
emits no edges is called a \emph{sink}. A \emph{path} $\mu$ in a graph $E$ is a sequence of edges
$\mu=e_1\dots e_n$ such that $r(e_i)=s(e_{i+1})$ for $i=1,\dots,n-1$. In this case, $s(\mu):=s(e_1)$ is the
\emph{source} of $\mu$, $r(\mu):=r(e_n)$ is the \emph{range} of $\mu$, and $n$ is the \emph{length} of $\mu$,
i.e, $l(\mu)=n$. We denote by $\mu^0$ the set of its vertices, that is:
$\mu^0=\{s(e_1),r(e_i):i=1,\dots,n\}$.

Now let $K$ be a field and let $KE$ denote the $K$-vector space which has as a basis the set of paths. It is
possible to define an algebra structure on $KE$ as follows: for any two paths $\mu=e_1\dots e_m, \nu=f_1\dots
f_n$, we define $\mu\nu$ as zero if $r(\mu)\neq s(\nu)$ and as $e_1\dots e_m f_1\dots f_n$ otherwise. This
$K$-algebra is called the \textit{path algebra} of $E$ over $K$.

For a field $K$ and a row-finite graph $E$, the {\em Leavitt path $K$-algebra} $L_K(E)$ is defined as the
$K$-algebra generated by a set $\{v\mid v\in E^0\}$ of pairwise orthogonal idempotents, together with a set
of variables $\{e,e^*\mid e\in E^1\}$, which satisfy the following relations:

(1) $s(e)e=er(e)=e$ for all $e\in E^1$.

(2) $r(e)e^*=e^*s(e)=e^*$ for all $e\in E^1$.

(3) $e^*e'=\delta _{e,e'}r(e)$ for all $e,e'\in E^1$.

(4) $v=\sum _{\{ e\in E^1\mid s(e)=v \}}ee^*$ for every $v\in E^0$ that emits edges.

Relations (3) and (4) are called of Cuntz-Krieger.

 The elements of $E^1$ are called \emph{(real) edges}, while for $e\in E^1$ we call $e^\ast$ a \emph{ghost
edge}.  The set $\{e^*\mid e\in E^1\}$ will be denoted by $(E^1)^*$.  We let $r(e^*)$ denote $s(e)$, and we
let $s(e^*)$ denote $r(e)$. If $\mu = e_1 \dots e_n$ is a path, then we denote by $\mu^*$ the element $e_n^*
\dots e_1^*$ of $L_K(E)$.

There exists a natural inclusion of the path algebra $KE$ into the Leavitt path algebra $L_K(E)$ sending
vertices to vertices and edges to edges. We will use this monomorphism without any explicit mention to it.

It is shown in \cite{AA1} that $L_K(E)$ is a ${\mathbb Z}$-graded $K$-algebra, spanned as a $K$-vector space
by $\{pq^* \mid p,q$ are paths in $E\}$. In particular, for each $n\in \mathbb{Z}$, the degree $n$ component
$L_KE)_n$ is spanned by elements of the form $pq^*$ where $l(p)-l(q)=n$.

The set of \emph{homogeneous elements} is $\bigcup_{n\in {\mathbb Z}} L_K(E)_n$, and an element of $L_K(E)_n$
is said to be $n$-\emph{homogeneous} or \emph{homogeneous of degree} $n$.

Note that the natural monomorphism from the path algebra $KE$ into the Leavitt path algebra $L_K(E)$ is
graded, hence $KE$ is a $\Z$-graded subalgebra of $L_K(E)$.

 An easy result which will be used later is the following one.

\begin{lem}\label{independencialineal} Any set of different paths
is $K$-linearly independent.
\end{lem}
\begin{proof}  Consider a graph $E$ and let $\mu_1, \dots, \mu_n$ be
different paths. Write $\sum_i k_i\mu_i=0$, for $k_i\in K$. Applying that $L_K(E)$ is $\Z$-graded we may
suppose that all the paths have the same length. Since $\mu_j^\ast\mu_i=\delta_{ij}r(\mu_j)$ then $0=\sum_i
k_i\mu_j^\ast\mu_i=k_jr(\mu_j)$; this implies $k_j=0$.
\end{proof}

\section{Algebras of right quotients}
 We first study the semiprimeness of the path algebra associated to a graph $E$.
Recall that an algebra $A$ is said to be \textit{semiprime} if it has no nonzero ideals of zero square,
equivalently, if $aAa=0$ for $a\in A$ implies $a=0$ (an algebra $A$ that satisfies this last condition is
called in the literature \emph{nondegenerate}).

\begin{prop}\label{SemipOfTheGraphAlgebra}
For a graph $E$ and a field $K$ the path algebra $KE$ is semiprime if and only if for every path $\mu$ there
exists a path $\mu^\prime$ such that $s(\mu^\prime)=r(\mu)$ and $r(\mu^\prime)=s(\mu)$.
\end{prop}
\begin{proof}
Suppose first that $KE$ is semiprime. Given a path $\mu$, since $\mu(KE)\mu\neq 0$, there exists a path
$\nu\in KE$ such that $\mu\nu\mu\neq 0$. This means that $s(\nu)=r(\mu)$ and $r(\nu)=s(\mu)$.

Now, let us prove the converse. Note that by \cite[Proposition II.1.4 (1)]{NvO}, a $\Z$-graded algebra is
semiprime if and only if it is graded semiprime. Hence, and taking into account that being graded semiprime
and graded nondegenerate are equivalent, it suffices to show that if $x$ is any nonzero homogeneous element
of $KE$, then $x(KE)x\neq 0$. Write $x=\sum_{i=1}^nk_i\alpha_i$, with $0\neq k_i\in K$ and $\alpha_1,\dots,
\alpha_n$ different paths of the same degree (i.e. of the same length). Denote the source and range of
$\alpha_1$ by $u_1$ and $v_1$, respectively. Then, by (3), ${\alpha_1}^\ast
x=k_1{\alpha_1}^\ast{\alpha_1}=k_1v_1$. By the hypothesis, there exists a path ${\alpha_1}^\prime$ such that
$s({\alpha_1}^\prime)=v_1$ and $r({\alpha_1}^\prime)=u_1$. Observe that ${\alpha_1}^\prime x\neq 0$;
otherwise $0={({\alpha_1}^\prime)}^\ast{\alpha_1}^\prime x=u_1x$, a contradiction since a set of different
paths is always linearly independent over $K$ (Lemma \ref{independencialineal}) and $\alpha_1=u_1\alpha_1\neq
0$. Therefore $0\neq k_1{\alpha_1}^\prime x=k_1v_1{\alpha_1}^\prime x = {\alpha_1}^\ast x{\alpha_1}^\prime
x\in {\alpha_1}^\ast x(KE)x$.
\end{proof}

If $a \in L_K(E)$ and $d \in {\mathbb Z}^+$, then we say that $a$ is \textit{representable as an element of
degree} $d$ \textit{in real (respectively ghost) edges} in case $a$ can be written as a sum of monomials from
the spanning set $\{pq^*\ \vert \ p, q \ \hbox{are paths in }E\}$, in such a way that $d$ is the maximum
length of a path $p$ (respectively $q$) which appears in such monomials. Note that an element of $L_K(E)$ may
be represented as an element of different degrees in real (respectively ghost) edges.

The $K$-linear extension of the assignment $pq^* \mapsto qp^*$ (for $p,q$ paths in $E$) yields an involution
on $L_K(E)$, which we denote simply as ${}^*$. Clearly $(L_K(E)_n)^* = L_K(E)_{-n}$ for all $n\in {\mathbb
Z}$.

Let $R \subseteq Q$ be rings. Recall that $Q$ is said to be an \textit{algebra of right quotients of} $R$ if
given $p, q\in Q$, with $p\neq 0$, there exists an element $r\in R$ such that $pr\neq 0$ and $qr\in R$.

\begin{prop}\label{rqa}
For any graph $E$ and any field $K$ the Leavitt path algebra $L_K(E)$ is an algebra of right quotients of the
path algebra $KE$.
\end{prop}
\begin{proof}
Consider $x, y \in L_K(E)$, with $x\neq 0$. Apply the reasoning in the first paragraph of the proof of
\cite[Proposition 3.1]{AMMS} to find a path $\mu$ such that $0\neq x\mu\in KE$. If $y\mu\in KE$ we have
finished the proof. If $y\mu\notin KE$, write $y\mu$ as a sum of monomials of the form $ke_1\dots
e_rf_1^\ast\dots f_s^\ast$ (it is possible by \cite[Lemma 1.5]{AA1}), with $k\in K$ and $e_1, \dots, e_r,
f_1, \dots, f_s \in E^1$, being this expression minimal in ghost edges, and denote by $m$ the maximum natural
number of ghost edges appearing in $y\mu$.

Note that $r(\mu)$ is not a sink; otherwise, for any $f\in E^1$,
$f^\ast r(\mu)=0$, hence $y\mu\in KE$, but this is not our case.
Therefore, there exists $h_1\in E^1$ such that $s(h_1)=r(\mu)$.
Then $x\mu h_1\neq 0$ because $x\mu\in KE\setminus\{0\}$ and by
Lemma \ref{independencialineal}. If $y\mu h_1\in KE$, our proof is
complete. Otherwise, repeat this reasoning. This process must stop
in at most $m$ steps.
\end{proof}

\section{Fountain-Gould and Moore-Penrose right orders}

Localization consists of assigning inverses to certain elements. This procedure can be carried out by taking
inverses (when the algebra has an identity) or,  more generally, by considering generalized inverses like
Moore-Penrose or group inverses, in which case the existence of the identity element plays no role at all.
There is a well-developed theory of rings of quotients when this kind of inverses is taken. We are speaking
about Fountain-Gould orders and Moore-Penrose orders.

For a graph $E$ having a not finite number of vertices, the Leavitt path algebra $L_K(E)$ is not unital,
hence it cannot be a maximal algebra of quotients of any of its subalgebras. However, we will show that it is
a Fountain-Gould right algebra of quotients of any of its semiprime subalgebras containing the path algebra
$KE$ when the graph $E$ is acyclic (a Moore-Penrose right order if we take into account the involution).

Moreover, we will describe in this section the structure of semiprime subalgebras of $L_K(E)$ containing the
path algebra $KE$, for $E$ an acyclic graph such that the set of line points is ``dense" in the sense that
the saturated closure of $P_l(E)$ is the whole set of vertices of $E$.

\begin{prop}\label{dmaximalcasofinito} Let $E$ be a finite and acyclic graph and let $A$ be a semiprime algebra
such that $KE \subseteq A\subseteq L_K(E)$. Then $L_K(E)=Q_{max}^r(KE)=Q_{cl}^r(A)$.
\end{prop}
\begin{proof} Apply \cite[Proposition 3.4]{E1} to the opposite graph of $E$
(the one obtained from $E$ by changing ranges to sources and sources to ranges) and \cite[Proposition
3.5]{AAS1} to show that $Q_{max}^r(KE)=L_K(E)$. The fact of $L_K(E)$ be semiprime and artinian and an algebra
of right quotients of $A$ (because it is an algebra of right quotients of $KE$, by Proposition \ref{rqa})
implies, by \cite[Corollary 3.4]{GS}, that it is the classical algebra of right quotients of $A$.
\end{proof}

Let $ a $ be an element in a ring $ R. $ We say that $ b $ in $ R  $ is the  \textit{group inverse} of $ a $
if the following conditions hold: $ aba=a, \quad bab=b,  \quad ab=ba. $

It is easy to see that $ a $ has a group inverse $ b $ in $R$ if and only if there exists a unique idempotent
$ e $ $ ( e=ab ) $ in $ R $ such that $ a $ is invertible in the  ring $ eRe $ (with inverse $ b $), hence
the group inverse is unique and $a$ is said to be \textit{locally invertible}. Denote by $a^\sharp$ the group
inverse of $a$.

An element $ a \in R $ is called \textit{square cancellable} if $ a^2 x = a^2 y  $ implies $ ax=ay $  and
$xa^2=ya^2$ implies $xa=ya$, for $ x, y \in R\cup \{ 1\}$ (for $x=1$ or $y=1$ this means that $a^2=a^2y$ or
$a^2x=a^2$ implies $a=ay$ or $ax=a$, respectively, and analogously for the right hand side). Denote by
$\mathcal{S}(R)$ the sets of all square cancellable elements of $R$.

 Recall that a subring $ R $ of a ring $Q$ is  a \textit{Fountain-Gould right order in} $ Q $ if:
\begin{enumerate}
\item Every element of $\mathcal{S}(R) $ has a group inverse in $ Q$ and
\item  every element $ q \in Q $ can be written in the form $ q=ab^\sharp ,$ where $ b \in \mathcal{S}(R) $ and
            $ a \in R.$
\end{enumerate}

Now we extend Proposition \ref{dmaximalcasofinito} to the non-necessarily finite case.

\begin{prop}\label{fountaingould}
Let $E$ be an acyclic graph and let $A$ be a semiprime algebra such that $KE \subseteq A\subseteq L_K(E)$.
Then $A$ is a Fountain-Gould right order in $L_K(E)$.
\end{prop}

\begin{proof}
Use \cite[Lemma 2.2]{AMP} to find a family $\{E_n\}$ of finite subgraphs of $E$ such that
$L_K(E)=\underrightarrow{\textrm{lim}}\ L_K(E_n)$. Since $E$ is acyclic, the same can be said about $E_n$. We
show first that if $a$ is a semiregular element of $A$, then it is locally invertible in $L_K(E)$. Indeed,
let $m$ be in $\mathbb{N}$ such that $a \in L_K(E_m)$. The element $a$ is semiregular in $L_K(E_m)$;
moreover, by \cite[Proposition 3.5]{AAS1}, $L_K(E_m)$ is a semisimple and artinian algebra, which implies
that $a$ is locally invertible in $L_K(E_m)$ (by \cite[Proposition 2.6]{FG1}), hence it is locally invertible
in $L_K(E)$.

Now, take $x\in L_K(E)$, and let $r$ be in $\mathbb{N}$ such that $x\in L_K(E_r)$. By Proposition
\ref{dmaximalcasofinito}, $L_K(E_r)$ is the classical algebra of right quotients of $A\cap L_K(E_r)$, so we
may write $x=ab^\sharp$, with $a, b\in A\cap L_K(E_r)$, where $b^\sharp$ denotes the inverse of $b$ in
$L_K(E_r)$. Then $b^\sharp$ is just the generalized inverse of $b$ in $L_K(E)$, and our claim has been
proved.
\end{proof}

\begin{rema}{\rm It is not possible to eliminate the hypothesis of semiprimeness in the previous results
 (Propositions \ref{dmaximalcasofinito} and \ref{fountaingould}).
Consider, for example, the following graph:
$$\xymatrix{{\bullet} \ar[r]  & {\bullet}}$$ Then $ \mathcal{M}_2(K)\cong L_K(E)=Q_{max}^r(KE)$ would imply,
applying \cite[Theorem 4.6]{GS}, that $KE$ is semiprime, which is false by virtue of Proposition
\ref{SemipOfTheGraphAlgebra}.}
\end{rema}

 There is another way of localizing when the ring has an involution.
In this case, it is possible to consider Moore-Penrose inverses instead of group inverses and $\ast$-square
cancellable elements instead of the square-cancellable ones. Rings of quotients with respect to Moore-Penrose
inverses were introduced in a recent paper of the author (\cite{S}), where it is shown among others that
these two ``local'' approaches are equivalent when the ring of quotients (in any of the two senses) is
semiprime and coincides with its socle. A general theory of this kind of localization is developed in
\cite{AM3}. Since the Leavitt path algebra has an involution $\ast$, we may consider $(A, \ast)$ as an
algebra with involution such that $KE\subseteq A\subseteq L_K(E)$, and then rewrite the previous result as
follows.

\begin{prop}\label{fountaingouldstar}
Let $E$ be an acyclic graph and let $(A, \ast)$ be a semiprime algebra with involution such that $KE
\subseteq A\subseteq L_K(E)$. Then $A$ is a Moore-Penrose right $\ast$-order in $L_K(E)$.
\end{prop}
\begin{proof} Use Proposition \ref{fountaingould}, the Main Theorem of \cite{S} and
\cite[Theorem 2.1]{AM3}.
\end{proof}

We define a relation $\ge$ on $E^0$ by setting $v\ge w$ if there is a path $\mu\in E^*$ with $s(\mu)=v$ and
$r(\mu)=w$. In this case we will say that $v$ \textit{connects to} the vertex $w$.

Some consequences can be derived from the results in \cite{GS} concerning subalgebras of $L_K(E)$ containing
the path algebra $KE$ for an acyclic graph $E$ such that every vertex connects to a line point. First, recall
the notion of line point (it appears in \cite{AMMS}).

A subset $H$ of $E^0$ is called \emph{hereditary} if $v\ge w$ and $v\in H$ imply $w\in H$. A hereditary set
is \emph{saturated} if every vertex which feeds into $H$ and only into $H$ is again in $H$, that is, if
$s^{-1}(v)\neq \emptyset$ and $r(s^{-1}(v))\subseteq H$ imply $v\in H$.

 The set $T(v)=\{w\in E^0\mid v\ge w\}$ is the \emph{tree} of $v$, and it is the smallest hereditary subset of
$E^0$ containing $v$. We extend this definition for an arbitrary set $X\subseteq E^0$ by $T(X)=\bigcup_{x\in
X} T(x)$. The \emph{hereditary saturated closure} of a set $X$ is defined as the smallest hereditary and
saturated subset of $E^0$ containing $X$. It is shown in \cite{AMP} that the hereditary saturated closure of
a set $X$ is $\overline{X}=\bigcup_{n=0}^\infty \Lambda_n(X)$, where
\begin{enumerate}
\item[] $\Lambda_0(X)=T(X)$, and
\item[] $\Lambda_n(X)=\{y\in E^0\mid
s^{-1}(y)\neq \emptyset$ and $r(s^{-1}(y))\subseteq \Lambda_{n-1}(X)\}\cup \Lambda_{n-1}(X)$, for $n\ge 1$.
\end{enumerate}

If $\mu$ is a path in a graph $E$ such that $s(\mu)=r(\mu)$ and $s(e_i)\neq s(e_j)$, for every $i\neq j$,
then $\mu$ is called a \emph{cycle}. We say that a vertex $v$ in $E^0$ is a \emph{bifurcation} (or that
\emph{there is a bifurcation at} $v$) if $s^{-1}(v)$ has at least two elements. A vertex $u$ in $E^0$ will be
called a \emph{line point} if there are neither bifurcations nor cycles at any vertex $w\in T(u)$, where
$T(v)=\{w\in E^0\mid v\ge w\}$ is the \emph{tree} of $v$. We will denote by $P_l(E)$ the set of all line
points in $E^0$.
 We say that a path $\mu$ \emph{contains no bifurcations} if the set
$\mu^0\setminus\{r(\mu)\}$ contains no bifurcations, that is, if none of the vertices of the path $\mu$,
except perhaps $r(\mu)$, is a bifurcation.

The notion of line point is central in the study of the socle of a Leavitt path algebra (see \cite{AMMS}):
for a graph $E$ the algebra $L_K(E)$ has a nonzero socle (equivalently, it has nonzero minimal left (right)
ideals) if and only if it has line points. In fact, the socle of $L_K(E)$ is generated, as an ideal, by the
line points in $E$.

Every ring $R$ which is left nonsingular and such that every element $a \in R$ has finite left Goldie
dimension (the left Goldie dimension of the left ideal generated by $a$ is finite) will be called a
\textit{left local Goldie} ring (equivalently, by \cite[5, (7.5)]{Lam}, $R$ satisfies the ascending chain
condition on the left annihilators of the form $\text{lan}(a)$, with $a \in R$, and every element $a \in R$
has finite left Goldie dimension). If additionally $R$ has finite left (global) dimension, then $R$ will be
called a \textit{left Goldie} ring.

\begin{theor}\label{leftlocalgoldie} Let $E$ be an acyclic graph
such that $\overline{P_l(E)}=E^0$. Then, for every semiprime algebra $A$ such that $KE\subseteq A\subseteq
L_K(E)$ we have:
\begin{enumerate}
\item[(i)] $A$ is a left local Goldie ring
\item[(ii)] $A$ is prime if and only it the only hereditary and
saturated subsets of $E^0$ are $\emptyset$ and $E^0$.
\item[(iii)] $A$ has finite left Goldie dimension if and only if
$E$ is finite.
\end{enumerate}
\end{theor}
\begin{proof}  (i).
The condition $\overline{P_l(E)}=E^0$ implies that $L_K(E)$ coincides with its socle (\cite[Theorem
 4.2]{AMMS}). Moreover, $L_K(E)$ is an algebra of right quotients of $A$ because
it is an algebra of right quotients of $KE$ (Proposition \ref{rqa}). Now, (i) follows by \cite[Proposition
3.5 (i)]{GS}. The statement (ii) can be obtained applying \cite[Proposition 3.5 (1)]{GS} and \cite[Theorem
3.11]{AA1}. Finally, (iii) follows from \cite[Proposition 3.5 (2)]{GS} and the obvious fact that $L_K(E)$
artinian implies $E$ finite.
\end{proof}

\section{Maximal algebras of quotients}

A Leavitt path algebra has essential socle if and only if every vertex of $E^0$ connects to a line point.
This characterization is the key tool to compute maximal algebras of right quotients of $L_K(E)$, for $E$
satisfying the above condition.

An edge $e$ is an {\it exit} for a path $\mu = e_1 \dots e_n$ if there exists $i$  such that $s(e)=s(e_i)$
and $e\neq e_i$.

\begin{prop}\label{nonsing}
For any graph $E$ the Leavitt path algebra $L_K(E)$ is nonsingular.
\end{prop}
\begin{proof} Suppose that the left singular ideal $Z_l(L_K(E))$ contains a nonzero element $x$. By
\cite[Proposition 3.1]{AMMS} there exist $\gamma, \mu \in L_K(E)$ such that $0 \neq \gamma x \mu \in Kv$ for
a certain vertex $v$ or $0 \neq \gamma x \mu \in wL_K(E)w \cong K[t, t^{-1}]$, where $w$ is a vertex for
which there is a cycle without exits based at it. Since the left singular ideal is an ideal, we have in the
first case that it contains the vertex $v$, which is not possible because $Z_l(L_K(E))$ does not have
idempotents. In the second case, $0 \neq \gamma x \mu \in Z_l(L_K(E))$ implies, using \cite[GS, proposition
2.1 (viii)]{GS} that the local algebra of $L_K(E)$ at $\gamma x \mu$, i.e. $L_K(E)_{\gamma x \mu}$, is left
nonsingular. But $L_K(E)_{\gamma x \mu}=(wL_K(E)w)_{\gamma x \mu}$, that is, there exists a nonzero element
$u\in K[t, t^{-1}]$ satisfying that $K[t, t^{-1}]_u$ is left nonsingular. This implies, by \cite[Proposition
2.1 (vii)]{GS}, that $u$ lies in the left singular ideal of $K[t, t^{-1}]$, a contradiction since $K[t,
t^{-1}]$ is a nonsingular algebra.

   The right nonsingularity of $L_K(E)$ can be proved analogously.
\end{proof}

\begin{rema} \label{density}
 {\rm The previous result implies that for any graph $E$ the essential left(respectively
right/bi-) modules over the Leavitt path algebra $L_K(E)$ coincide with the dense left(respectively
right/bi-) modules (use \cite[Lemma 2.24 (b)]{G}). Similar results can be obtained if we consider the graded
notions.}
\end{rema}

\begin{theor}\label{socleessential}
Let $E$ be a graph. Then $Soc(L_K(E))$ is an essential ideal of $L_K(E)$ if and only if every vertex connects
to a line point.
\end{theor}
\begin{proof} Suppose first that every vertex connects to a line point. Let $x$ be a nonzero element in $L_K(E)$.
By \cite[Proposition 3.1]{AMMS} there exist $v\in E^0$ and $\gamma, \mu\in L_K(E)$ such that $0\neq \gamma x
\mu=kv\in Kv$, or there exists a cycle $c$ without exits, and $w\in c^0$, such that $0\neq \gamma x \mu\in
wL_K(E)w\cong K[t, t^{-1}]$. In the first case, since every vertex connects to a line point, there exist
$u\in P_l(E)$ and a path $\alpha\in E^\ast$ satisfying $s(\alpha)=v$ and $r(\alpha)=u$. Then
$u=\alpha^\ast\alpha=\alpha^\ast v\alpha=k^{-1}\alpha^\ast\gamma x \mu\alpha$. The second case cannot happen
because since $w$ connects to an element of $P_l(E)$ and the vertices of any cycle are not in $P_l(E)$ (by
the very definition), $w$ connects to a vertex which is not in $c^0$, hence $c$ has an exit.  This shows that
$Soc(L_K(E))$ is an essential ideal of $L_K(E)$.

Now, suppose that the socle is an essential ideal of $L_K(E)$. We remark that every vertex of
$\overline{P_l(E)}$ connects to a line point because $\overline{P_l(E)}= \bigcup_{n=0}^\infty
\Lambda_n(P_l(E))$ and since $P_l(E)$ is a hereditary subset of $E^0$, every $\Lambda_n(P_l(E))$ is
hereditary again. Pick a vertex $v$ in $E^0$. Apply that $Soc(L_K(E))$ is essential as a left ideal, and
\cite[Theorem 4.2]{AMMS}, to find $\alpha, \alpha_i, \beta_i\in L_K(E)$ and $u_i\in P_l(E)$, for $i=1, \dots,
r$, such that $0\neq \alpha v=\sum_i \alpha_iu_i\beta_i=\sum_i \alpha_iu_i\beta_i v$, hence for some $i$ we
have $0\neq u_i\beta_i v$. Suppose for a moment that only (1) can happen in \cite[Proposition 3.1]{AMMS}.
Then, there exist: an element $\mu \in KE$, $k\in K$ and $w\in E^0$ such that $0\neq u_i\beta_i v\mu=kw$.
Note that $w\in \overline{P_l(E)}$ because $ u_i \in Soc(L_K(E))=I(\overline{P_l(E)})$ (\cite[Theorem
4.2]{AMMS}) and by \cite[Lemma 2.1]{APS}. If $\gamma$ is a path connecting $w$ to a line point $u$, then
$0\neq \gamma =w\gamma u= k^{-1} u_i\beta_i v\mu\gamma u$. In particular, $v\mu\gamma u\neq 0$. Take into
account that $\mu\gamma\in KE$, hence for some path $\lambda$ appearing in the summands of $\mu\gamma$,
$0\neq v\lambda u$, that is, $\lambda$ is a path that joins $v$ to $u$.

Finally, justify that every cycle has an exit.  If $c$ were a cycle without exits, and $v$ any vertex in
$c^0$, with similar ideas to that of \cite[Proof of Theorem 3.11]{AA1} it could be shown that $vL_K(E)v\cong
K[x, x^{-1}]$. Apply that the socle of $L_K(E)$ is essential as a left ideal to find a nonzero $\alpha\in
L_K(E)v\cap Soc(L_K(E))$. By the semiprimeness of $L_K(E)$ (\cite[Proposition 1.1]{AMMS}) there exists $\beta
\in L_K(E)$ satisfying $0\neq y:= v\beta\alpha v\in vL_K(E)v\cap Soc(L_K(E))$. Then, since $K[x, x^{-1}]$ has
no minimal one-sided ideals the local algebra  of $vL_K(E)v$ at $y$, that is, $(vL_K(E)v)_y=L_K(E)_y$, is not
artinian by \cite[Proposition 2.1 (v)]{GS} (see \cite{GS} for the definition of local algebra at an element).
But, on the other hand, the same result and $y\in Soc(L_K(E))$ imply that $L_K(E)_y$ is artinian, a
contradiction.
\end{proof}

The following definitions are particular cases  of those appearing in \cite[Definition 1.3]{DHSz}:

Let $E$ be a graph, and let $\emptyset \ne H\in \mathcal{H}_E$. Define
$$F_E(H)=\{ \alpha =\alpha_1 \dots \alpha_n  \mid \alpha _i\in
E^1, s(\alpha _1)\in E^0\setminus H, r(\alpha _i)\in E^0\setminus H \mbox{ for } i<n, r(\alpha _n)\in H\}.$$
Denote by $\overline{F}_E(H)$ another copy of $F_E(H)$. For $\alpha\in F_E(H)$, we write $\overline{\alpha}$
to denote a copy of $\alpha$ in $\overline{F}_E(H)$. Then, we define the graph ${}_HE=({}_HE^0, {}_HE^1, s',
r')$ as follows:
\begin{enumerate}
\item $({}_HE)^0=H\cup F_E(H)$. \item $({}_HE)^1=\{ e\in
E^1\mid s(e)\in H\}\cup \overline{F}_E(H)$.
\item For every $e\in E^1
\mbox{ with } s(e)\in H$, $s'(e)=s(e)$ and $r'(e)=r(e)$. \item For every $\overline{\alpha}\in
\overline{F}_E(H)$, $s'(\overline{\alpha})=\alpha$ and $r'(\overline{\alpha})=r(\alpha)$.
\end{enumerate}

A (semiprime) algebra which coincides with its socle will be called \textit{semisimple}.

\begin{theor}\label{elmaximaleseldelzocalo}
Let $E$ be a graph such that every vertex connects to a line point and denote by $H$ the saturated closure of
$P_l(E)$. Then:
 $$Q_\upsilon(L_K(E))\cong Q_\upsilon(I(H))\cong Q_\upsilon(L_K(_HE)),$$
  where $Q_\upsilon(-)$ is
the Martindale symmetric algebra of quotients, the maximal symmetric algebra of quotients, the maximal
algebra of left/right quotients or the graded maximal algebra of left/right quotients, and $_HE$ is an
acyclic graph such that $\overline{P_l(_HE)}=_HE^0$, that is, $L_K(_HE)$ is a semisimple algebra.
\end{theor}
\begin{proof}
Note that by \cite[Theorem 4.2]{AMMS}, $Soc(L_K(E))=I(H)$ (where $I(H)$ is the ideal of $L_K(E)$ generated by
the vertices of $H$; it is a graded ideal of $L_K(E)$ -see \cite[Remark 2.2]{APS}-), and by \cite[Lemma
1.2]{AP}, $I(H)\cong L_K(_HE)$, where $_HE$ is an acyclic graph (see the proof of \cite[Theorem 4.6]{AMMS})
and satisfies that the saturated closure of $P_l(_HE)$ is $(_HE)^0$ ($L_K(_HE)$ coincides with its socle and
we apply \cite[Theorem 4.2]{AMMS}). By Theorem \ref{socleessential} and Remark \ref{density}$I(H)$ is a dense
subalgebra of $L_K(E)$, hence the result follows for the Martindale symmetric, the maximal symmetric and the
maximal one-sided algebras of quotients. For the graded case, apply \cite[Lemma 2.8]{ASgraded}.
\end{proof}

\begin{rema}\label{ejemplografo}
{\rm The condition ``every vertex connects to a line point'' is milder than
$\overline{P_l(E)}=E^0$. For example, the graph
$$
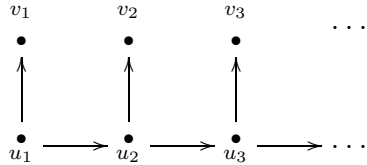
\xymatrix{
\underset{{\bullet}}{^{v_1}} & \underset{{\bullet}}{^{v_2}} & \underset{{\bullet}}{^{v_3}}  &{\dots}{}& \\
\ar[u]\overset{{\bullet}}{_{u_1}} \ar[r] & \ar[u]\overset{{\bullet}}{_{u_2}}  \ar[r]& \ar[r]
\ar[u]\overset{{\bullet}}{_{u_3}}& \ar@{.}\dots}$$ satisfies the first condition, although not the second
one, as was shown in \cite[Example 4.6]{AMMS}.}
\end{rema}

\begin{corol}\label{elmaximaldelalgebradecaminos}
Let $E$ be a graph such that every vertex connects to a line point and denote by $H$ the saturated closure of
$P_l(E)$. Then
$$Q_{max}^r(KE)= Q_{max}^r(L_K(E))\cong Q_{max}^r(I(H))\cong Q_{max}^r(L_K(_HE)),$$
 where $_HE$ is
an acyclic graph such that $\overline{P_l(E)}=E^0$, that is, $L_K(_HE)$ is a semisimple algebra.

The same can be said about the maximal graded algebra of right quotients.
\end{corol}
\begin{proof} Apply Theorem \ref{elmaximaleseldelzocalo} and
Proposition \ref{rqa}.
\end{proof}

\begin{exem} {\rm Consider the graph $E$ in Remark \ref{ejemplografo}. Since every vertex
connects to a line point, by Theorem \ref{socleessential}, $I(H)=Soc(L_K(E))$ is an essential ideal of
$L_K(E)$. It was shown in \cite[Example 4.6]{AMMS} that $Soc(L_K(E))=I(H)$, where $H=\{v_n \ \vert\ n \in
\mathbb{N}\}=P_l(E)$. Apply Theorem \ref{elmaximaleseldelzocalo} to have that $Q_\upsilon(L_K(E))\cong
Q_\upsilon(L_K(_HE))$.  Define $e_n$ as the element of $E^1$ having source $u_n$ and range $v_n$, and by
$f_n$ the edge whose source is $u_n$ and whose range is $u_{n+1}$, for every $n \in \mathbb{N}$. Then $_HE$
is the graph:}

$$
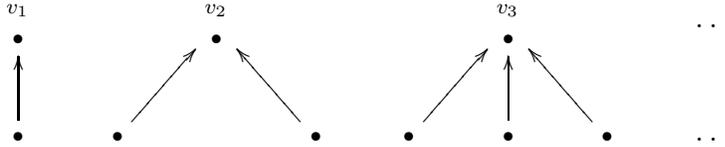
\xymatrix{
\underset{{\bullet}}{^{v_1}} &  &\underset{{\bullet}}{^{v_2}}& & &\underset{{\bullet}}{^{v_3}}&  &\dots& \\
\ar[u]\overset{{\bullet}}  {}& \ar[ur]\overset{{\bullet}}{}& &\ar[ul]\overset{{\bullet}}{}  &
\ar[ur]\overset{{\bullet}}{}& \ar[u]\overset{{\bullet}}{}&\ar[ul] \overset{{\bullet}}{}& \ar@{.}\dots}$$

 {\rm It is clear that the associated Leavitt path algebra is isomorphic to $\underset{n=2}{\overset{\infty}\oplus} M_n(K)$,
therefore $Q_{\upsilon}(L_K(E))\cong Q_\upsilon(\underset{n=2}{\overset{\infty}\oplus}
M_n(K))=\underset{n=2}{\overset{\infty}\prod} M_n(K)$.}
\end{exem}

\section{Toeplitz algebras}

The Toeplitz algebra is defined (see, for example \cite{HS}) as
the C$^*$-algebra of continuous functions on the quantum disc. It
can be described as the graph algebra associated to the following
graph (see, for example, \cite{HMS}):

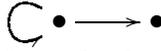
$$\xymatrix{ {\bullet} \ar@(ul,dl) \ar[r]&{\bullet}}$$
\medskip
\noindent although other descriptions in terms of graph algebras can be given (see \cite[Section 3]{HS}).

\begin{defi}\label{toeplitzalgebra} {\rm Define the (\textit{algebraic}) \textit{Toeplitz algebra} as the
Leavitt path algebra associated to the graph given above. Denote
it by $T$.}
\end{defi}

Note that the Toeplitz algebra (we will avoid the use of the word
``algebraic'' because we are in an algebraic context) is a Leavitt
path algebra associated to a graph $E$ for which every vertex
connects to a line point. Our main concern in this section will be
to give a description of the Toeplitz algebra, similar to that
given in the analytic context, and to show that if we substitute
the line point in the graph for a connected acyclic graph without
bifurcations, and the edge connecting the loop to the sink to any
(finite) number of edges, then the resulting Leavitt path algebra
is again the defined Toeplitz algebra.

We first start by describing those Leavitt path algebras $L_K(E)$ such that there are neither bifurcations
nor cycles at any point of $E^0$.

As we have said before, it was shown in \cite[Lemma 1.5]{AA1} that every monomial in $L_K(E)$ is of the form:
$kv$, with $k\in K$ and $v\in E^0$, or $k e_1\dots e_m f_1^*\dots f_n^*$ for $k\in K$, $m, n\in \mathbb{N}$,
$e_i,f_j\in E^1$. By a \textit{reduced expression of a monomial} we will understand an expression of the form
$e_1\dots e_m f_1^*\dots f_n^*$ with $m+n$ minimal. If this is the case, we will say that $e_1\dots e_m
f_1^*\dots f_n^*$ \textit{is a reduced monomial}.

The following definition can be found in \cite[pg. 56]{R}: a \emph{walk} in a directed graph $E$ is a path in
the underlying undirected graph. Formally, a walk $\mu$ is a sequence $\mu=\mu_1\dots \mu_n$ with $\mu_i\in
E^1\cup (E^1)^*$ and $s(\mu_i)=r(\mu_{i+1})$ for $1\leq i < n$. The directed graph $E$ is \emph{connected} if
for every two vertices $v,w\in E^0$ there is a walk $\mu=\mu_1\dots \mu_n$ with $v=s(\mu)$ and $w=r(\mu)$.
Intuitively, $E$ is connected if $E$ cannot be written as the union of two disjoint subgraphs, or
equivalently, $E$ is connected in case the corresponding undirected graph of $E$ is so in the usual sense. It
is easy to show that if $E$ is the disjoint union of subgraphs  $\{E_i \}$, then $L_K(E)\cong \oplus
L_K(E_i)$. If all the $E_i$'s are connected, then each $E_i$ will be called a \textit{connected component} of
$E$.

\begin{prop}\label{basiswithnobifurcations} Let $E$ be an acyclic graph such that there are no bifurcations at any
vertex of $E^0$ and write $E=\bigcup_{i\in \Upsilon}E_i$, where the $E_i$'s are the connected components of
$E$. Then $L_K(E)$ is isomorphic to $\bigoplus_{i\in\Upsilon}\mathcal{M}_{\alpha_i}(K)$, where
$\alpha_i=\sharp(E_i^0)$, the cardinal of $E_i^0$, and  $\alpha_i\in \mathbb{N} \cup\{\infty\}$.
\end{prop}
\begin{proof}
As we have explained, $L_K(E)\cong \bigoplus_{i\in\Upsilon}L_K(E_i)$, hence we may reduce our study to the
case of a connected graph $E$. We have divided the proof into three steps.

{\underline{Step 1.}} Every monomial in $L_K(E)$ has a unique reduced expression, that is, for every monomial
$z\in L_K(E)$, there exist two unique paths $\alpha, \beta\in L_K(E)$, with $length(\alpha)+lenght(\beta)$
minimal such that $z=\alpha\beta^*$.

 Let $\alpha, \beta, \mu, \nu$ be paths such that $\alpha\beta^*=\mu\nu^*$, being both reduced
expressions. Suppose that the length of $\alpha$ is strictly less to the length of $\mu$. Use (3) to show
that if $\alpha=e_1 \dots e_n$, then $\mu=e_1 \dots e_ng_1\dots g_r$, for some $e_i, g_j\in E^1$, and write
$\beta^*=f_1^*\dots f_m^*$ and $\nu^*=h_1^*\dots h_{m+r}^*$, with $f_i, h_j\in E^1$. Then $\quad e_1 \dots
e_n f_1^*\dots f_m^*=e_1 \dots e_ng_1\dots g_rh_1^*\dots h_{m+r}^*$ implies, by (3), $f_1^*\dots
f_m^*=g_1\dots g_rh_1^*\dots h_{m+r}^*$. Since there are no bifurcations in $E$, $(f_m\dots f_1)(f_1^*\dots
f_m^*)=s(f_m)$. Using this and (3), after multiplying the expression below by $f_m\dots f_1$ on the left hand
side, we obtain:

$$f_m\dots f_1 g_1\dots g_r h_1^*\dots h_{m+r}^*=s(f_m).$$

   This means that $s(f_m)=r(h_{m+r}^*)=s(h_{m+r})$, but there are no bifurcations in $E$, hence
$f_m=h_{m+r}$. Multiply by $f_m^*$ on the left hand side and by $f_m$ on the right hand side. Then

$$f_{m-1}\dots f_1 g_1\dots g_r h_1^*\dots h_{m+r-1}^*=s(f_{m-1}).$$

Proceed again in this form and in $m+r$ steps we will have obtained $f_i=h_{r+i}$, for $i\in \{1, \dots, m
\}$ and $g_j=h_{r-j+1}$ for $j\in\{1, \dots, r\}$. This implies that

 $$\mu\nu^*=e_1 \dots e_ng_1\dots g_rh_1^*\dots h_{m+r}^*= e_1 \dots e_nh_r\dots h_1h_1^*\dots h_{m+r}^*=
 e_1 \dots e_nh_{r+1}^*\dots h_{m+r}^*,$$
\noindent
 which is not a reduced expression of $\mu\nu^*$ since here the number of edges plus the number of ghost edges is
$n+m$, while in the first expression of $\mu\nu^*$ this sum was $n+m+2r$.

Consequently, $r=0$ and $e_1 \dots e_n f_1^*\dots f_m^*=e_1 \dots e_nh_1^*\dots h_m^*$. This implies (by (4)
and since there are no bifurcations) $f_1^*\dots f_m^*=h_1^*\dots h_m^*$; by (3), $f_i=h_i$ for $i\in \{1,
\dots,m\}$, and so $\alpha=\mu$ and $\beta=\nu$.

{\underline{Step 2.}} The set of all reduced monomials is a basis of $L_K(E)$ as a $K$-vector space. Denote
it by $\mathcal{B}$.

 We know from~\cite[Lemma 1.5]{AA1} that the expressions $\{\alpha\beta^*\}$ generate $L_K(E)$ as a vector
space, hence we only need to check that they are linearly independent.

Suppose we have $\sum_i k_i\alpha_i\beta_i^*=0$, where all the summands are different from zero and
$\alpha_i\beta_i^*\neq\alpha_j\beta_j^*$. Moreover, taking into account the degree, we may suppose that each
summand has the same degree.

 Let $\beta_1$ be with maximal length among the $\beta_i$'s. Then $0=\sum_i k_i\alpha_i\beta_i^*\beta_1$,
 where at least one of the summands is nonzero ($k_1\alpha_1\beta_1^*\beta_1=k_1\alpha_1$) and being each
 summand in only real edges. Let us call $\gamma_1=\alpha_1$ and $\gamma_i=\alpha_i\beta_i^*\beta_1$, for $i\neq
 1$. With this notation the formula below reads $\sum_i \gamma_i=0$. Multiply by $\gamma_1^*$ on the right
 hand side and apply that $\gamma_1\gamma_1^*=s(\gamma_1)$, because there are no exits. Then multiply by
 $\gamma_1^*$ on the  left hand side and apply (3). We obtain the following formulas:

 $$-k_1s(\gamma_1)=\sum_{i\neq 1}k_i\gamma_i\gamma_1^*$$

 $$-k_1r(\gamma_1)=\sum_{i\neq 1}k_i\gamma_1^*\gamma_i$$

The first one implies that $s(\gamma_i)=s(\gamma_1)$ for all $i$'s given nonzero terms. The second one that
$r(\gamma_i)=s(\gamma_1)$ for the same i's. This implies $\gamma_i=\gamma_1$ for these terms. Hence, there
exists at least one $i$ such that $\alpha_1=\gamma_1=\gamma_i=\alpha_i\beta_i^*\beta_1$, so
$\alpha_1\beta_1^*=\alpha_i\beta_i^*$, a contradiction.

{\underline{Step 3.}} The result.

Let $u$ and $v$ be vertices. Since we are considering that the graph $E$ is connected, there exists a walk
$\mu$ starting at $v$ and ending at $w$. By Step 1, $\mu$ has a reduced expression and it is an element of
the basis $\mathcal{B}$ (Step 2), therefore we may describe $\mathcal{B}$ as the set of walks $\mu_{j,k}$,
where $\mu_{jk}$ is the only element in $L_K(E)$ such that $s(\mu_{jk})=v_j$ and $r(\mu_k)=v_k$, for $v_j,
v_k\in E^0$. Define, for $\alpha=\sharp(E^0)$ the map $\varphi: L_K(E)\to \mathcal{M}_\alpha(K)$ as the
$K$-linear map which acts on the elements of $\mathcal{B}$ as follows: $\varphi(\mu_{jk})=e_{jk}$, being
$e_{jk}$ the matrix unit having all the entries equal to zero except those in row $j$ and column $k$. It is
not difficult to show that it is an isomorphism of $K$-algebras.
\end{proof}

\begin{theor}\label{descriptionoftoeplitzalgebras}
Let $F$ be an acyclic graph such that every vertex connects to a
line point, and consider the graph $E=E(n, F)$ such that
$E^0=F^0\cup \{v\}$, $E^1=F^1\cup\{e, e_1, \dots, e_n\}$, where
$e$ has range and source $v$, $s(e_i)=v$, $r(e_i)\in F^0$, and for
every $f\in F^0$, its range and source are the corresponding as
edges of the graph $F$. Then:
\begin{enumerate}
\item[(i)] $Soc(L_K(E))$ is an essential ideal of $L_K(E)$.
\item[(ii)] There exists an exact sequence:
$$0 \rightarrow Soc(L_K(E))\rightarrow T\rightarrow K[x, x^{-1}]\rightarrow 0.$$
\item[(iii)] $L_K(E)$ is isomorphic to the Toeplitz algebra $T$.
\item[(iv)] There is a subalgebra $T^\prime$ of $\mathrm{RCFM}(K)$,
isomorphic to $T$, such that
$$\mathcal{M}_\infty(K)\subseteq T^\prime \subseteq \mathrm{RCFM}(K).$$
\end{enumerate}
\end{theor}
\begin{proof} (i). Since every vertex of $E^0$ connects to a line
point, we may apply Theorem \ref{socleessential} to obtain that the socle of $L_K(E)$ is an essential ideal.

(ii). Observe that $F^0$, which is the saturated closure of the set of line points of $E$, is an hereditary
and saturated subset of $E^0$, therefore $Soc(L_K(E))=I(F^0)$ by \cite[Theorem 5.2 ]{AMMS}. Moreover,
\cite[Lemma 2.3 (1)]{APS} implies that
$$0 \rightarrow Soc(L_K(E))\rightarrow T\rightarrow L_K(E/F^0)\rightarrow 0$$
is an exact sequence, where $E/F^0$ is the graph having one vertex and one edge ($E/F^0$ is called the
quotient graph; see, for example, \cite{AMP} for its definition). Since $L_K(E/F^0)\cong K[x, x^{-1}]$, we
have proved this item.

(iii). By \cite[Lemma 1.2]{AP}, $I(F^0)$ is isomorphic to $L_K(_{F^0}E)$. Note that $_{F^0}E$ is an acyclic
graph without bifurcations at any point and having an infinite number of edges. In the particular case of
being $n=1$ and $F$
 the graph having one vertex and no edges, $_{F^0}E$ is
$$
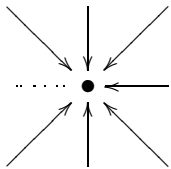
\xymatrix{ \ar[dr]& \ar[d]& \ar[dl]\\& \ar@{.} @(l,l)\bullet&\ar[l] \\ \ar[ur]& \ar[u]& \ar[ul]}$$

 By Proposition \ref{basiswithnobifurcations},
$L_K(_{F^0}E)$ is isomorphic to $\mathcal{M}_\infty(K)$. This fact, jointly with (ii), imply the result.

(iv). It has been proved in (iii) that $I(F^0)\cong L_K(_{F^0}E)\cong \mathcal{M}_\infty(K)$. By Corollary
\ref{elmaximaldelalgebradecaminos} $Q_{max}^r(L_K(E))\cong Q_{max}^r(I(F^0))\cong
Q_{max}^r(\mathcal{M}_\infty)= \mathrm{RCFM}(K).$
\end{proof}

%%%%%%%%%%%%%%%%%%%%%%%%%%%%%%%%%%%%%%%%%%%%%%%%%%%%%%%%%%%%%%%%%%%%%%%%%%%%%%%%%%%%%%%%%%%%%%%%%%%%%%%%%%%%%%%
%%%%%%%%%%%%%%%%%%%%%%%%%%%%%%%%%%%%%%%%%%%%%%%%%%%%%%%%%%%%%%%%%%%%%%%%%%%%%%%%%%%%%%%%%%%%%%%%%%%%%%%%%%%%%%%
\section*{acknowledgments}
Partially supported by the Spanish MEC and Fondos FEDER through project MTM2004-06580-C02-02 and by the Junta
de Andaluc\'{\i}a and Fondos FEDER, jointly, through project FQM-336.
%%%%%%%%%%%%%%%%%%%%%%%%%%%%%%%%%%%%%%%%%%%%%%%%%%%%%%%%%%%%%%%%%%%%%%%%%%%%%%%%%%%%%%%%%%%%%%%%%%%%%%%%%%%%%%%
%%%%%%%%%%%%%%%%%%%%%%%%%%%%%%%%%%%%%%%%%%%%%%%%%%%%%%%%%%%%%%%%%%%%%%%%%%%%%%%%%%%%%%%%%%%%%%%%%%%%%%%%%%%%%%%

\end{document}